\documentclass{amsart}

\theoremstyle{remark}
\newtheorem{teo}{\it Theorem}
\newtheorem{define}{\it Definition}
\newtheorem{ejem}{\it Example}
\newtheorem{lema}{\it Lemma}
\newtheorem{cor}{\it Corollary}
\newtheorem{prop}{\it Proposition}

\begin{document}

\newcommand{\slq}{\ensuremath{(sl_{n+1}^{+})_{q}\;}}
\newcommand{\A}{\ensuremath{\mathcal{A}}}
\newcommand{\wL}{\ensuremath{\widetilde{L}}}
\newcommand{\R}{\ensuremath{\mathcal{R}}}
\newcommand{\B}{\ensuremath{\mathcal{B}}}
\newcommand{\eref}[1]{\eqref{#1}}

\title[ Quantum Lie algebras, PBW and the Yang-Baxter equation]{Quantum Lie algebras of type $\mathbb{A}_{n}$ positive, PBW bases
and the Yang-Baxter equation}
\author[C. Bautista]{C\'esar Bautista}
\address{Facultad de Ciencias de la Computaci\'on,
Benem\'erita Universidad Aut\'onoma de Puebla, Edif. 135. 14 sur y Av.
San Claudio, Ciudad Universitaria, Puebla, Puebla, C.P.
72570, M\'exico}

\email{bautista@fcfm.buap.mx, bautista@solarium.cs.buap.mx}

\thanks{This work was partially supported by a grant of Fundaci\'on TELMEX
and by the
Investigation Project INI106897  of DGAPA/UNAM}

\begin{abstract}
We show explicitly a generalised Lie algebra embedded in the positive and
negative parts of the Drinfeld-Jimbo quantum groups of type
$\mathbb{A}_{n}$. Such a generalised Lie algebra satisfy axioms closely
related to the ones found by S.L. Woronowicz.
For the universal enveloping algebra of such generalised Lie algebras we
establish several conditions in order to obtain bases of type
Poincar\'e-Birkhoff-Witt. Besides a graded algebra is
proposed and some relations with the quantum Yang-Baxter equation are studied.
\end{abstract}

%\pacs{02.20.Sv, 02.10.Vr, 02.10.Tq}

%\submitted

%Available at Los Alamos e-Print server: {\tt http://xxx.lanl.gov}

\maketitle

\begin{section}{Introduction}
The applications of the theory of Lie algebras and groups to
theoretical physics are well known. However theoretical problems have
reached a level which requires more general algebraic objects. For
instance, the objects called {\it quantum groups} \cite{Chaichian} which
contains Lie groups and Lie algebras as degenerate cases. Moreover, related
to quantum groups appear different kinds of generalised Lie algebras
\cite{Woronowicz, Majid, Vybornov}. In particular, inside
of the Drinfeld-Jimbo quantum groups (or quantised universal
enveloping algebras) several authors \cite{Lyubashenko, Delius, yo} have
discovered generalised Lie algebras, or {\it quantum Lie algebras}. In this
paper we deal with the Drinfeld-Jimbo quantum group $U_{q}(sl_{n+1})$.

In \cite{yo} is shown that there exists inside  $U_{q}(sl_{n+1})$ a
generalised color Lie algebra \slq such that
its universal enveloping algebra is the positive part
$U_{q}^{+}(sl_{n+1})$ of  $U_{q}(sl_{n+1})$.
Unfortunately the axioms of such a generalised
Lie algebra depend on the chosen basis.

The aim of this paper is to propose generalised Lie algebra axioms for
the module \slq independent of the choice of basis. We call such a
structure {\it $T$-Lie algebra}. A part of such axioms
(antisymmetry, Jacobi identity) are
closely related to the ones found by Woronowicz \cite{Woronowicz}.

Besides we shall prove that \slq  satisfies an additional Jacobi identity in a
similar manner to the {\it balanced generalised Lie algebras} of
Lyubashenko-Sudbery \cite{Lyubashenko}. This means, \slq is
a balanced generalised Lie algebra. However, the universal enveloping
algebra of \slq as a $T$-Lie algebra is not the same universal
enveloping algebra as a balanced generalised Lie algebra.

Finally the relations with the quantum Lie algebra definition of
Vybornov \cite{Vybornov} are studied. Following the work of Vybornov  we
obtain that, in some sense, the quantum Yang-Baxter equation (QYBE) includes
not only a generalised Jacobi identity, and a generalised antisymmetry, but
an additional property called multiplicativity,
\pagebreak
which is a part of the
conditions in order to obtain bases of type Poincar\'e-Birkhoff-Witt (PBW).
\end{section}

\begin{section}{Basic Notions and Common Constructions}
Let $k$ be a commutative unitary ring.
\begin{define}\label{defT}
A $T$-Lie algebra is a $k$-module $L$ with morphisms
\begin{align*}
[,]:L\otimes_{k}L&\rightarrow L,\quad\mbox{bracket}\\
S:L\otimes_{k}L&\rightarrow L,\quad\mbox{presymmetry}\\
\langle,\rangle:L\otimes_{k}L&\rightarrow
L\otimes_{k}L,\quad\mbox{pseudobracket}
\end{align*}
such that the following axioms are satisfied
\begin{enumerate}
\item{Stability}
   \begin{enumerate}
   \item There exist $L=\oplus_{\eta\in\mathbb{N}}L_{\eta}$ a strict
    gradation of $L$ related to $[,]$, (this
    means $[L_{\eta_{1}},L_{\eta_{2}}]\subseteq
    L_{\eta_{1}+\eta_{2}-1},\;\forall\,\eta_{1},\eta_{2}$).
    \item \[\langle L_{\eta_{1}},L_{\eta_{2}}\rangle\subseteq
    (L\otimes_{k}L)_{\eta_{1}+\eta_{2}-1},\quad\forall\eta_{1},\eta_{2}
           \]
       where
       $(L\otimes_{k}L)_{\eta}=\sum_{\eta_{1}+\eta_{2}\leq\eta}L_{\eta_{1}}
       \otimes_{k}L_{\eta_{2}}$.
    \end{enumerate}
\item{Antisymmetry}
      \begin{enumerate}
      \item $[,]S=-[,]$
      \item $\langle,\rangle S=-\langle,\rangle$
      \item $[,]\langle,\rangle=0$
      \end{enumerate}
\item{Jacobi identity}\label{jacobi}
\begin{equation}\label{jT}
[,](1\otimes [,])-[,]([,]\otimes 1)-[,](1\otimes [,])(S\otimes
1)=0
\end{equation}
\end{enumerate}
\end{define}
\begin{ejem}
A $G$-algebra $\Lambda$ is a graded $k$-module
$\Lambda=\oplus_{n}\Lambda_{n}$ equipped with two multiplications,
$(\lambda,\gamma)\mapsto \lambda\gamma$  and $(\lambda,\gamma)\mapsto
[\lambda,\gamma]$ with some additional  properties, we are interested
just in the properties of the bracket $[,]$. These are
\begin{equation}
[\lambda,\gamma]\in \Lambda_{m+n-1},\mbox{ if
}\lambda\in\Lambda_{m},\gamma\in\Lambda_{n}\label{st}
\end{equation}
\begin{equation}
[\lambda,\gamma]=-(-1)^{(m-1)(n-1)}[\gamma,\lambda],\mbox{ if
}\lambda\in\Lambda_{m},\gamma\in\Lambda_{n},\label{ant}
\end{equation}
and if
$\lambda\in\Lambda_{m},\gamma\in\Lambda_{n},\eta\in\Lambda_{p}$,
\begin{equation}
(-1)^{(m-1)(p-1)}[\lambda,[\gamma,\eta]]
+(-1)^{(n-1)(m-1)}[\gamma,[\eta,\lambda]]
+(-1)^{(p-1)(n-1)}[\eta,[\lambda,\gamma]]=0\label{jG}
\end{equation}
The property \eref{st} says that $\Lambda=\oplus_{n}\Lambda_{n}$ is a
strictly graded algebra related to $[,]$. Therefore $\Lambda$ has a
structure of $T$-Lie algebra with bracket $[,]$, pseudobracket
$\langle,\rangle=0,$  and presymmetry $S$ defined by
\[
S(\lambda\otimes\gamma)=(-1)^{(m-1)(n-1)}\gamma\otimes\lambda,\mbox{ if
}\lambda\in\Lambda_{m},\gamma\in\Lambda_{n}.
\]
Since,
\begin{eqnarray*}
S(1\otimes [,])=(1\otimes [,])(1\otimes S)(S\otimes 1)\\
S([,]\otimes 1)=([,]\otimes 1)(S\otimes 1)(1\otimes S)
\end{eqnarray*}
then the Jacobi identity \eref{jG} is equivalent to the Jacobi identity
\eref{jT} of $\Lambda$ as a $T$-Lie algebra.

Observe that we are not using the standard {\it reduced degree} of
$\Lambda$ (see \cite{Gerstenhaber}, p 89).
\end{ejem}
If $L$ is a $T$-Lie algebra we define its {\it universal enveloping
algebra} $U(L)$ as the quotient
\[
U(L)=L^{\otimes}\,/J
\]
where $L^{\otimes}$ is the $k$-tensorial algebra of $L$ and $J$ is the
two sided ideal  generated by
\[
x\otimes y-S(x\otimes y)-\langle x,y\rangle-[x,y],\quad\forall x,y\in L.
\]

Related to each $T$-Lie algebra $L$ there exist an {\it abelian $T$-Lie
algebra} denoted $L^{0}$. Its structure is: $L^{0}=L$ as $k$-module,
$[,]^{0}=0,$ $\langle,\rangle^{0}=0,$ $S^{0}=S$. The universal enveloping
algebra
\[
\mathcal{S}(L)=U(L^{0})
\]
is called {\it $q$-symmetric algebra of $L$}.

Now, let $\mathcal{B}$ a totally ordered basis of $L$ a $T$-Lie algebra
with strict gradation given by
\[
L=\oplus_{\eta\in\mathbb{N}}L_{\eta}.
\]
Let $x_{\lambda}\in\mathcal{B},\Sigma=(x_{\lambda_{1}},\ldots,x_{\lambda_{u}})$
 finite-non-decreasing sequence of elements of $\mathcal{B}$. We write
 $\eta(\lambda)=\eta(x_{\lambda})=\alpha$ if $x_{\lambda}\in L_{\alpha}$,
 $z_{\lambda}=x_{\lambda}\in\mathcal{S}(L),$
 $z_{\Sigma}=z_{\lambda_{1}}\ldots z_{\lambda_{u}}\in\mathcal{S}(L),$
 $z_{\emptyset}=1\in\mathcal{S}(L)$,
 $\eta(\Sigma)=\eta(z_{\Sigma})=\eta(x_{\lambda_{1}})+\ldots
 +\eta(x_{\lambda_{u}})$. Besides we put $x_{\lambda}\leq \Sigma\;$ if
 $x_{\lambda}\leq x_{\lambda_{1}}$.

 Let $\mathcal{S}_{p}$ be the $k$-submodule of
$\mathcal{S}(L)$ generated by $z_{\Sigma}$ such that $\eta(\Sigma)\leq
p$.

\begin{lema}\label{acc}
Let $L$ be a $T$-Lie algebra and $\mathcal{B}$ a free basis of $L$ as
$k$-module, $\mathcal{B}$ totally ordered, such that
\begin{enumerate}
\item $S(x\otimes y)=q_{x,y}y\otimes x,\;q_{x,y}\in
k,\forall x,y\in \mathcal{B}$;
\item $q_{x,y}=q^{-1}_{y,x}$.
\end{enumerate}
Then, there exist a $k$-morphism
\[
\underline{\quad}\cdot\underline{\quad}:L\otimes_{k}\mathcal{S}(L)\rightarrow
\mathcal{S}(L)
\]
satisfying
\begin{enumerate}
\item[(A)] $x_{\lambda}\cdot z_{\Sigma}=z_{\lambda}z_{\Sigma}$ for
$x_{\lambda}\leq \Sigma;$
\item[(B)] $x_{\lambda}\cdot z_{\Sigma}-z_{\lambda}z_{\Sigma}\in
\mathcal{S}_{\eta(\lambda)+\eta(\Sigma)-1}$.
\end{enumerate}
\end{lema}
\begin{proof}
See \cite{yo}, proof of lemma V.1.
\end{proof}
\end{section}
\begin{section}{Poincar\'e-Birkhoff-Witt bases and $T$-Lie algebras}
In \cite{yo} we proved that the PBW theorem is not
a general property for all the basic $T$-Lie algebras. We have to
restrict them. Such a restriction is called {\it adequate}. Thanks to
lemma \ref{acc} the condition {\it adequate} makes sense not only for
{\it basic} $T$-Lie algebras but for $T$-Lie algebras.

\begin{define}
An adequate $T$-Lie algebra is a $T$-Lie algebra  with a
totally ordered free basis for $L$ as a $k$-module such that they hold the
condition {\it adequate} (see \cite{yo}, definition V.1).
\end{define}

We do not have to confuse {\it adequate basic $T$-Lie algebra} (defined
in \cite{yo}) with
{\it adequate $T$-Lie algebra} (defined in this paper).

Let $L^{3}$ be the $k$-submodule of the $k$-tensorial algebra
$L^{\otimes}$ generated by
\[
x_{i}\otimes x_{j}\otimes x_{k},\;x_{i}<x_{j}<x_{k}\in\mathcal{B}
\]
and let $\;^{3}L$ be the $k$-submodule generated by
\[
x_{k}\otimes x_{j}\otimes x_{i},\;x_{i}<x_{j}<x_{k}\in\mathcal{B}.
\]
Let us put $S_{2},S_{1}:L^{3\otimes}\rightarrow L^{3\otimes}$,
$S_{1}=1\otimes S$, $S_{2}=(1\otimes S)$.
\begin{prop}
Let $L$ be an adequate $T$-Lie algebra related to the basis
$\mathcal{B}$ such that,
\begin{enumerate}
\item \label{m1}$S([,]\otimes 1)|_{L^{3}}=(1\otimes[,])S_{1}S_{2}|_{L^{3}};$
\item \label{m2}$S(1\otimes [,])|_{L^{3}}=([,]\otimes 1)S_{2}S_{1}|_{L^{3}};$
\item \label{m3}$S([,]\otimes 1)|_{\;^{3}L}=(1\otimes [,])S_{1}S_{2}|_{\;^{3}L};$
\item \label{m4}$S(1\otimes [,])|_{\;^{3}L}=([,]\otimes 1)S_{2}S_{1}|_{\;^{3}L};$
\end{enumerate}
then $U(L)$ holds a basis of type PBW related to $\mathcal{B}$.
\end{prop}
\begin{proof}
According with \cite{yo}, it will suffice to prove that $L$ is a basic $T$-Lie
algebra. The conditions \eref{m1}, \eref{m2} ensure multiplicativity, while
conditions \eref{m3}, \eref{m4} say that Jacobi identity for $L$ as a
basic $T$-Lie algebra
follows from the Jacobi identity as $T$-Lie algebra. Therefore $L$ is a
basic $T$-Lie algebra.
\end{proof}
\begin{prop}\label{prejacobi}
Let $L$ be a basic $T$-Lie algebra with basis $\mathcal{B}$. If it holds
\begin{enumerate}
\item \label{c1}$S([,]\otimes 1)|_{\;^{3}L}=(1\otimes [,])S_{1}S_{2}|_{\;^{3}L};$
\item \label{c2}$S(1\otimes [,])|_{\;^{3}L}=([,]\otimes 1)S_{2}S_{1}|_{\;^{3}L};$
\item \label{ad}$[,]S([x_{i},x_{j}]\otimes x_{j})=[,](1\otimes
[,])S_{1}S_{2}(x_{i}\otimes x_{j}\otimes x_{j}),$ $\forall
x_{i},x_{j}\in\mathcal{B}.$
\end{enumerate}
and a second Jacobi identity
\begin{equation}\label{jacobi2}
[,](\,(1\otimes [,])S_{1}S_{2}- ([,]\otimes 1)S_{2}S_{1}+([,]\otimes
1)S_{2}\,)|_{\;^{3}L}=0
\end{equation}
then
\begin{eqnarray}
\label{first}[,](1\otimes [,])-[,]([,]\otimes 1)-[,](1\otimes
[,])S_{1}=0\\
\label{second}[,](1\otimes [,])-[,]([,]\otimes 1)+[,]([,]\otimes
1)S_{2}=0.
\end{eqnarray}
\end{prop}

\begin{proof}
Let $x,y,z\in\mathcal{B}$. If $x>y>z$ then the Jacobi identity of $L$ as a
  $T$-Lie algebra together with the hypothesis \eref{c1} and \eref{c2} say
  \begin{equation}\label{j1}
  [x,[y,z]]-[[x,y],z]-q_{x,y}[y,[x,z]]=0
  \end{equation}
  besides, the second Jacobi identity \eref{jacobi2} can be written as
  \begin{equation}\label{j2}
  [x,[y,z]]-[[x,y],z]+q_{y,z}[[x,z],y]=0.
  \end{equation}

First, we are going to prove
\begin{equation}\label{j}
[x,[y,z]]-[[x,y],z]-q_{x,y}[y,[x,z]]=0.
\end{equation}
There are several cases:
\begin{enumerate}
\item[\underline{Case}] $x=y=z$: then \eref{j}  trivially holds.
\item[\underline{Case}] Two of $x,y,z$ are the same: now there are several
subcases;
  \begin{enumerate}
  \item[\underline{Subcase}] $x=y,z$: the left side of \eref{j} is
  $[x,[x,z]]-q_{x,x}[x,[x,z]]=0$.
  \item[\underline{Subcase}] $x,y=z$: the left side of \eref{j} is
  \begin{equation*}
  -[[x,y],y]-q_{x,y}[y,[x,y]]=[S([x,y]\otimes y]-q_{x,y}[y,[x,y]]=0
  \end{equation*}
  since hypothesis \ref{ad}.
  \item[\underline{Subcase}] $x=z,y$: again, the left side of \eref{j} is
  \begin{align*}
  [x,[y,x]]-[[x,y],x]&=[x,[y,x]]+q_{x,y}[[y,x],x]\\
                     &=[x,[y,x]]-q_{y,x}[\,S([y,x]\otimes x)\,]\\
                     &=[x,[y,x]]-q_{y,x}q_{x,y}[x,[y,x]]
  \end{align*}
  because hypothesis \eref{ad}.
  \end{enumerate}
\item[\underline{Case}] $x,y,z$ three different elements: let be $s$ the
left side of \eref{j}. There are several subcases:
$x>y>z,x>z>y,y>x>z,y>z>x,z>x>y,z>y>x$
  \begin{enumerate}
  \item[\underline{Subcase}] $x>y>z$: equation \eref{j1}
  \item[\underline{Subcase}] $x>z>y$:
  \begin{align*}
  s&=-q_{y,z}[x,[z,y]]-[[x,y],z]+q_{x,y}[y,[x,z]]\\
  &=-q_{y,z}[x,[z,y]]-[[x,y],z]
  -q_{x,y}q_{y,x}q_{y,z}[[x,z],y]\\
  &=-q_{y,z}([x,[z,y]]-[[x,z],y]+q_{z,y}[[x,y],z])=0
  \end{align*}
thanks to \eref{j2}.
  \item[\underline{Subcase}] $y>x>z$:
  \begin{align*}
  s&=[x,[y,z]]+q_{x,y}[[y,z],z]-q_{x,y}[y,[x,z]]\\
   &=-q_{x,y}([y,[x,z]]-[[y,x],z]-q_{y,x}[x,[y,z]])=0
  \end{align*}
  since the Jacobi identity in $L$.
  \item[ ] The remaining subcases are similar.
  \end{enumerate}
\end{enumerate}

By similar computing \eref{second} holds too.
\end{proof}
\end{section}
\begin{section}{A quantum Lie algebra of type $\mathbb{A}_{n}$ positive}
Let $[,]$ be the usual bracket of the Lie algebra $sl_{n+1}$,
$\mathcal{B}=\{\,e_{ij}\,|\, 1\leq i<j\leq
n+1\,\}$ canonical basis of the Lie subalgebra $sl_{n+1}^{+}$ form by upper
triangular matrices. Put $h_{ij}=[e_{ij},e_{ij}^{t}],$
$1\leq i<j\leq n+1$. Then $[h_{ab},e_{ij}]=c_{ab}e_{ij},$ where
$c_{ab,ij}\in\mathbb{Z}$.
Let $k$ be an unitary commutative ring,  \slq  the $k$-free module with basis
$\mathcal{B}$. Define: an
order by
\[
e_{ab}<e_{uv},\mbox{ if }a+b<i+j\mbox{ or }(a+b=i+j\mbox{ and }b<j),
\]
a $k$-morphism $S:\slq\otimes_{k}\slq\rightarrow\slq\otimes_{k}\slq$ by
\[
S(e_{ab}\otimes e_{ij})=q^{c_{ab,ij}}e_{ij}\otimes e_{ab},\mbox{ if
}e_{ab}<e_{ij},\;S(e_{ab}\otimes e_{ab})=e_{ab}\otimes e_{ab},\quad
S^{2}=1,
\]
a $k$-morphism  $[,]_{q}:\slq\otimes\slq\rightarrow\slq$ by
\[
[e_{ab},e_{ij}]_{q}=[e_{ab},e_{ij}]\mbox{ if
}e_{ab}<e_{ij},\quad[,]_{q}S=-[,]_{q}
\]
a $k$-morphism $\langle,\rangle:
\slq\otimes_{k}\slq\rightarrow\slq\otimes_{k}\slq$  by
\begin{equation*}
\mbox{ if }e_{ab}<e_{ij},\;\langle e_{ab},e_{ij}\rangle=\begin{cases}
                             (q-q^{-1})e_{ib}\otimes e_{aj}& \mbox{ if
                             $a<i<b$ and $i<b<j$,}\\
                             0 & \mbox{ otherwise,}
                             \end{cases}
\end{equation*}
and
\begin{equation*}
\langle,\rangle S=-\langle,\rangle.
\end{equation*}
\begin{lema}\label{lema}
If $e_{ab}>e_{ij}>e_{uv}$ then
\[
[[e_{uv},e_{ab}],e_{ij}]=[[e_{uv},e_{ab}]_{q},e_{ij}]_{q}=0.
\]
\end{lema}
\begin{proof}
Since $[,]_{q}$ is a deformation of the usual bracket $[,]$ of $sl_{n}$,
this amounts to showing that $[[e_{uv},e_{ab}],e_{ij}]=0$.

Because $e_{uv}<e_{ab}$ we can suppose $v=a$ then we have to prove
\begin{equation}\label{eq}
[e_{ub},e_{ij}]=0.
\end{equation}
Notice that if $b\neq i$ and $j\neq u$ then
\eref{eq} holds. But if $b=i$ then $a<b=i<j$ implies $a+i<j+i$ and
$e_{ab}=e_{ai}<e_{ij}$, a contradiction. And if $j=u$  then
$e_{ua}<e_{iu}$ this implies $u+a<i+u$ or $u+a=i+u$ and $a<u$, it
follows $a<u$, again a contradiction because $u<v=a$. Therefore
\eref{eq} holds.
\end{proof}

\begin{teo}
The $k$-module \slq has a structure of $T$-Lie algebra.
\end{teo}
\begin{proof}
We have to prove the Jacobi identity \eref{jacobi} of the definition
\ref{defT} for the bracket $[,]_{q}$. Since proposition \ref{prejacobi}, it
suffices to prove that its conditions hold. Let $e_{ij}$, $1\leq
i<j\leq n$ be canonical basis of \slq. If $e_{ab}<e_{ij}<e_{uv}$ then
\[
[e_{ab},e_{ij}]_{q}<e_{uv},\mbox{ if }[e_{ab},e_{ij}]_{q}\neq 0\mbox{
and }e_{ab}<[e_{ij},e_{uv}]_{q},\mbox{ if }[e_{ij},e_{uv}]_{q}\neq 0
\]
besides, if $[,]$ denotes the classical bracket of $sl_{n}$ and
$h_{ij}=[e_{ij},e_{ij}^{t}]$ then $[h_{ij},e_{uv}]=[h_{uv},e_{ij}]=
c_{ij,uv}e_{uv}$ for certain $c_{ij,uv}\in k$, and from the Jacobi
identity in $sl_{n}$ it follows
\[
[h_{ab},[e_{ij},e_{uv}]]=(c_{ab,ij}+c_{ab,uv})[e_{ij},e_{uv}]
\]
this implies the conditions \eref{c1} and \eref{c2} of proposition
\ref{prejacobi}.
The condition \eref{ad} of proposition \ref{prejacobi} holds because
$[[e_{ij},e_{ab}]_{q},e_{ab}]_{q}=0$, $\forall e_{ij},e_{ab}$. Now only
remains to prove the second Jacobi identity \eref{jacobi2}.

Let $x,y,z$ basic elements such that $x>y>z$, then the left side of
\eref{j2} is,
\begin{align*}
&=-q_{xy}q_{xz}q_{yz}([z,[y,x]]-[[z,y],x])-q_{yz}q_{xz}[[z,x],y]\\
&=q_{xz}q_{yz}(q_{x,y}[[z,x],y]-[[z,x],y])\\
&=0
\end{align*}
since lemma \ref{lema}.
\end{proof}
\begin{prop}
The $T$-Lie algebra \slq is a balanced generalised Lie algebra.
\end{prop}
\begin{proof}
Define $\gamma:(sl_{n}^{+})_{q}\otimes_{k} (sl_{n}^{+})_{q}\rightarrow
(sl_{n}^{+})_{q}\otimes_{k} (sl_{n}^{+})_{q}$ by $\gamma=1-S$. Then,
equations \eref{first} and \eref{second} of proposition \ref{prejacobi}
are the
left Jacobi identity and the right Jacobi identity respectively, (see
\cite{Lyubashenko}).
\end{proof}
\end{section}
\begin{section}{On the canonical graded algebra related to the universal enveloping
algebra}
For the $T$-Lie algebras we can repeat the classical construction \cite{Humphreys} (or color
\cite{Scheunert}) of the graded algebra related
to a universal enveloping algebra.

Let $L=\oplus_{\eta}L_{\eta}$ strictly gradation of the $T$-Lie algebra
$L$.  Define the $k$-submodules of $L^{\otimes}$ the tensorial algebra of $L$,
\begin{equation*}
T_{m}=\coprod_{\mbox{\tiny $\begin{array}{c}
               k\geq 0\\
              \eta_{1}+\ldots+\eta_{k}\leq
              m
              \end{array}$}}(L_{\eta_{1}}\otimes\ldots\otimes
              L_{\eta_{k}}),\;m\geq
              0
\end{equation*}
\begin{equation*}
T^{m}=\coprod_{\mbox{\tiny $\begin{array}{c}
               k\geq 0\\
              \eta_{1}+\ldots+\eta_{k}=
              m
              \end{array}$}}(L_{\eta_{1}}\otimes\ldots\otimes
              L_{\eta_{k}}),\;m\geq
              0
\end{equation*}
and let $\pi:L^{\otimes}\rightarrow U(L)$ natural morphism. Put
$U_{m}=\pi(T_{m})$, $U_{-1}=0$. Set the $k$-module
$G^{m}=U_{m}/U_{m-1}$, and let the multiplication on $U(L)$ define a
bilinear map $G^{m}\times G^{p}\rightarrow G^{m+p}$. This extends at
once to a bilinear map $G\times G\rightarrow G$ making $G=\coprod_{m\geq
0}G^{m}$ a graded
associative algebra with $1$.

Since $\pi$ maps $T^{m}$ into $U_{m}$, the composite linear map
$\varphi_{m}:T^{m}\rightarrow U_{m}\rightarrow G^{m}$ makes sense. It is
surjective because from $\pi(T_{m}-T_{m-1})=U_{m}-U_{m-1}$. Using the
universal property of the $k$-module coproduct
$L^{\otimes}=\coprod_{m\geq 0}T^{m}$ we get a $k$-morphism
$\varphi:L^{\otimes}\rightarrow G$, which is surjective.

Let $I$ be the $k$-bilateral ideal of $L^{\otimes}$ generated by
\[
x\otimes y-S(x\otimes y),\;\forall x,y\in L.
\]

\begin{prop} The map $\varphi:L^{\otimes}\rightarrow G$ is a $k$-algebra
morphism. Moreover, $\varphi(I)=0$, so $\varphi$ induces a $k$-algebra
epimorphism $\omega:\mathcal{S}(L)\rightarrow G$.
\end{prop}
\begin{proof}
First, notice that $\varphi:L^{\otimes}\rightarrow G$ is a $k$-algebra
morphism because the product definition of $G$.

Let $x\otimes y-S(x\otimes y),\; (x,y\in L)$ be a generator of
$I$ such that $x\in L_{\eta_{1}},y\in L_{\eta_{2}}$. Then
$\pi(x\otimes y-S(x\otimes y))\in U_{\eta}$ where
$\eta=\eta_{1}+\eta_{2}$. On the other hand $\pi(x\otimes y-S(x\otimes
y))=\pi([x,y]+\langle x,y\rangle)\in U_{\eta-1}$, whence
$\pi(x\otimes y-S(x\otimes y))\in U_{\eta-1}/U_{\eta-1}=0.$ It follows
$I\subseteq Ker\,\varphi$.
\end{proof}

In the classical case $\omega$ is an isomorphism. However, this is not
true in general, because if we take the $T$-Lie algebra
$L=\widetilde{(sl_{4})}_{q}$, (see \cite{yo}, examples III.5, IV.4 and
definition IV.1) then
$x_{2}x_{6}=0$ in $U(L)$. Besides if
$p:L^{\otimes}\rightarrow \mathcal{S}(L)$ is the natural projection then
$\omega(\,p(x_{2}\otimes x_{4})\,)=\pi(x_{2}\otimes x_{6})=0$.
Therefore $\omega$  it is not injective.

Although if $L$ is a basic $T$-Lie algebra, the symmetric algebra
$\mathcal{S}(L)$
is isomorphic to the graded algebra $G$.
\begin{cor}
\noindent
\begin{enumerate}
\item If $L$ is a basic $T$-Lie algebra then $\mathcal{S}(L)\simeq G$ as algebras.
\item If $L$ is a basic $T$-Lie algebra then $U(L)$ has no zero divisors
$\neq 0$.
\item The algebras $U\slq$ and $M_{p,q,\epsilon}(n)$ which is a
multiparametric deformation of $GL(n)$, have no zero divisors $\neq 0$.
\end{enumerate}
\end{cor}
\begin{proof}
\noindent
\begin{enumerate}
\item Just copy word by word the classical proofs \cite{Humphreys}, p
94 or
\cite{Jacobson}, p 166.
\item See \cite{Jacobson}, theorem 4, p 164.
\item The algebras $U\slq$ and $M_{p,q,\epsilon}(n)$ are both universal
enveloping algebras of basic $T$-Lie algebras, (see \cite{yo}).
\end{enumerate}
\end{proof}
\end{section}
\begin{section}{The Yang-Baxter equation and $T$-Lie algebras}
Now we are going to apply some remarks of Vybornov \cite{Vybornov} to the
theory of $T$-Lie algebras in order to obtain solutions of
the QYBE.

Let $L$ be a free $k$-module with basis $\B$ totally ordered and
$S:L\otimes_{k} L\rightarrow L\otimes_{k}L$ such that $S(x\otimes y)=
q_{x,y}y\otimes x$ for all $x,y\in\B$ and certain $q_{x,y}\in k$. Let us
put $\wL=L\oplus k$. We may extend $S$ to \wL\; by means of defining
\[
S(x\otimes 1)=1\otimes x,\;S(1\otimes
x)=x\otimes 1,\;S(1\otimes 1)= 1\otimes 1,\;\forall x\in L, 1\in k.
\]

Define two families of linear maps
\begin{gather*}
\mathcal{R}(\lambda),\mathcal{R}(\lambda)^{\prime}:\wL\otimes\wL\rightarrow
\wL\otimes \wL\\
\mathcal{R}(\lambda)=S+[,](p\otimes p)\otimes\lambda,\;
\mathcal{R}(\lambda)^{\prime}=S+\lambda\otimes [,](p\otimes p).
\end{gather*}
where $p:\wL\rightarrow L$ is the natural projection.

As usual, we define $\mathcal{R}_{1}(\lambda)=\mathcal{R}(\lambda)\otimes 1$,
etc.
\begin{prop}
Let $V$ be a $k$-submodule of $\wL^{3\otimes}$.
\begin{enumerate}
\item
$\R(\lambda)\R(\lambda)^{\prime}=\R(\lambda)^{\prime}\R({\lambda})=1,$
for any $\lambda$ if and only if $[,]S=-[,]$;
\item 
$\R_{1}(\lambda)\R_{2}(\lambda)\R_{1}(\lambda)|_{V}=
 \R_{2}(\lambda)\R_{1}(\lambda)\R_{2}(\lambda)|_{V}$, for any $\lambda$ if
 and only if the Jacobi identity \eref{jT} holds on $V$ and the following
 multiplicativity conditions hold also
 \begin{gather*}
 S(1\otimes [,])|_{V}=([,]\otimes 1)S_{1}S_{2}|_{V}\quad
 S([,]\otimes 1)|_{V}=(1\otimes [,])S_{2}S_{1}|_{V}\\
 S(1\otimes [,])S_{1}|_{V}=([,]\otimes 1)S_{2}|_{V}.
 \end{gather*}
 \item
 $\R_{1}(\lambda)^{\prime}\R_{2}(\lambda)^{\prime}
 \R_{1}(\lambda)^{\prime}|_{V}=
 \R_{2}(\lambda)^{\prime}\R_{1}(\lambda)^{\prime}
 \R_{2}(\lambda)^{\prime}|_{V}$, for any $\lambda$ if and only if the
 Jacobi identity \eref{second} holds on $V$ and the following multiplicativity
 conditions hold also
 \begin{gather*}
 S(1\otimes [,])|_{V}=([,]\otimes 1)S_{1}S_{2}|_{V}\quad
 S([,]\otimes 1)|_{V}=(1\otimes [,])S_{2}S_{1}|_{V}\\
 S([,]\otimes 1)S_{2}|_{V}=(1\otimes [,])S_{1}|_{V}.
 \end{gather*}
 \end{enumerate}
\end{prop}
\begin{proof}
By straightforward computations on the basic elements.
\end{proof}

\begin{cor}
Let $L$ be $\slq$ with canonical basis \B.
\begin{enumerate}
\item Let us put $x<0,\;\forall x\in\B$. Let $V_{1}$ be the $k$-submodule of
$L^{3}$ generated by
$x\otimes y\otimes z$ such that
$x,y,z\in\B$, $x<y<z$ and $y<[x,y]_{q}$. Then, for any $\lambda$,
\[
 \R_{1}(\lambda)\R_{2}(\lambda)\R_{1}(\lambda)|_{V_{1}}=
 \R_{2}(\lambda)\R_{1}(\lambda)\R_{2}(\lambda)|_{V_{1}}.
\]
\item Let us put $x>0,\;\forall x\in\B$. Let $V_{2}$ be the $k$-submodule of
$L^{3}$ generated by
$x\otimes y\otimes z$ such that
$x,y,z\in\B$, $x<y<z$ and $y>[x,y]_{q}$. Then, for any $\lambda$,
\[
\R_{1}(\lambda)^{\prime}\R_{2}(\lambda)^{\prime}
 \R_{1}(\lambda)^{\prime}|_{V_{2}}=
 \R_{2}(\lambda)^{\prime}\R_{1}(\lambda)^{\prime}
 \R_{2}(\lambda)^{\prime}|_{V_{2}}.
\]
\end{enumerate}
\end{cor}

Since the QYBE is equivalent to the braid equation   (see \cite{Larson}
p 3316),
we have obtained restricted solutions to the QYBE.

Besides, for $L=(sl_{3}^{+})_{q}$, we have
$V_{1}=V_{2}=L^{3\otimes}$. So, $(sl_{3}^{+})_{q}$ is a quantum Lie
algebra according with Vybornov \cite{Vybornov}.
\end{section}

\end{document}